\documentclass[12 pt, reqno]{amsart}
\usepackage{amsmath, amscd}

\title[]{Derived Functors Related to Wall Crossing}

\author[]{Kevin J. Carlin}

\address{Department of Mathematics and Computer Science\\
Assumption College\\
500 Salisbury St.\\
Worcester MA 01609-1296}

\email{kcarlin@assumption.edu}

\date{}


\renewcommand{\l}{\lambda}
\newcommand{\T}{\Theta}

\DeclareMathOperator{\Coker}{Coker}
\DeclareMathOperator{\Ker}{Ker}

\DeclareMathOperator{\im}{Im}
\DeclareMathOperator{\rad}{rad}
\DeclareMathOperator{\soc}{soc}
\DeclareMathOperator{\p}{pr}
\DeclareMathOperator{\h}{Hom}
\DeclareMathOperator{\E}{Ext}

\newcommand{\op}[2]{\mathop{#1_#2}}
\newcommand{\cop}[1]{\mathop{{}_{#1}\T}}
\newcommand{\kop}[1]{\mathop{\strut^{#1}\T}}
\newcommand{\dop}[1]{\mathop{\T\strut^{#1}}}

\newtheoremstyle{plainnopunct}%
{6pt plus 1pt minus 1pt}%
{6pt plus 1pt minus 1pt}%
{\itshape}%
{}%
{\bfseries}%
{}%
{.5em}%
{}

\theoremstyle{plainnopunct}
\newtheorem{thm}[equation]{Theorem}
\newtheorem{lem}[equation]{Lemma}
\newtheorem{prop}[equation]{Proposition}
\newtheorem{cor}[equation]{Corollary}
\newtheorem{conj}[equation]{Conjecture}
\newtheorem{ex}[equation]{Example}
\newtheorem{define}[equation]{Definition}

\calclayout
\advance\topmargin -\headheight
\advance\topmargin -\headsep
\advance\textheight \headheight
\advance\textheight \headsep
\numberwithin{equation}{section}

\begin{document}

\begin{abstract}
The setting is the representation theory of a simply connected, semisimple algebraic group over a field of positive characteristic. There is a natural transformation from the wall-crossing functor to the identity functor. The kernel of this transformation is a left exact functor. This functor and its first derived functor are evaluated on the global sections of a line bundle on the flag variety. It is conjectured that the derived functors of order greater than one annihilate the global sections. Also, the principal indecomposable modules for the Frobenius subgroups are shown to be acyclic.
\end{abstract}

\maketitle

\setcounter{section}{-1}

\section{Introduction}
\label{intro}

Let $G$ be a simply connected, semisimple algebraic group over an algebraically closed field $k$ of characteristic $p>0$. $B$ is a Borel subgroup containing a maximal torus $T$. $G_r$ denotes the kernel of the $r^\text{th}$ power of the Frobenius endomorphism. The character group of $T$ is $X(T)$. Assume that $p\ge h$ where $h$ is the Coxeter number.

Choose a dominant, $p$-regular weight $\l\in X(T)$. Let $s$ represent one of the walls of the alcove containing $\l$. The translation functor from $\l$ to the $s$-wall is denoted by $\op Ts$. Its adjoint functor is denoted by $\op Rs$. The wall-crossing functor $\op\T s$ is the composition $\op Rs{}\op Ts{}$. 

The co-unit $\epsilon$ of the adjunction $(\op Rs,\op Ts)$ is a natural transformation  from $\op\T s$ to the identity. The cokernel of $\epsilon$ is denoted by $\cop s{}$. The properties of this functor were the subject of a previous article \cite{KC2}. The functor of principal interest here is $\kop s{}=\Ker\epsilon$ and its derived functors. This functor has a left adjoint defined by $\dop s{}=\Coker\eta'$ where $\eta'$ is the unit of the adjunction $(\op Ts ,\op Rs)$ (proposition \ref{prop: adjpair}). 

The functors $\kop s{}$, $\dop s{}$, and $\cop s{}$ can be defined in a similar way in category $\mathcal O$ and its variants. Their properties and the calculation of their derived functors first appeared in \cite{KC1}. The functor $\dop s{}$ later appeared as a shuffle functor in Irving's work on shuffled Verma modules \cite[\S 3]{I}. More recently, Mazorchuk and Stroppel \cite[\S 5]{MS} described $\kop s{}$ as a co-shuffle functor. They also re-discovered the results of \cite{KC1} on the derived functors of $\kop s{}$ and $\dop s{}$. 

Based on the results in category $\mathcal O$ \cite[2.9ii]{KC1}, one might expect that the following properties should be true for all $G$ modules $V$ in the block of $\l$.
\begin{align}
\label{prop1}
R^1[\kop s{}]V&\simeq \cop sV\\
\label{prop2}
R\strut^q[\kop s{}]V&=0 \quad\text{ if }q>1
\end{align}
In fact, neither property holds true in general (see the remarks following \ref{ex: Lex}).

Using methods developed in \cite{KC2}, the derived functors of $\kop s{}$ can be determined in an indirect way (proposition \ref{prop: centralseq}). Under the condition that $p>2h-2$, this is used to show that the principal indecomposable $G_r$ modules are acyclic for $\kop s{}$ (corollary \ref{cor: Qracycle}).

The main result (proposition \ref{prop: R1}) shows that property \ref{prop1} is true if $V$ is the global sections of a line bundle over $G/B$ (under the same restriction on $p$). Using this result, the values of the derived functors $R\strut^q[\kop s{}]$ are determined for $q\le1$ (theorems \ref{thm: sRqI}, \ref{thm: R1s}, and \ref{thm: R0s}).  It is conjectured that the second property is also true for these modules (\ref{conj: der2vanish}).

The natural application of these derived functors is the calculation of extensions of $G$ modules. It is shown that some of the consequences of conjecture \ref{conj: der2vanish} are, at least, consistent with known results on extensions (proposition \ref{prop: consequences}). Finally, the possible existence of a spectral sequence corresponding to the composition $\h_G(V,\kop s(-))$ where $\eta'_V$ is an injective map is considered.

\section{Notations}
\label{nt}

The notations are close to the standard ones (as defined in \cite{J}). The exceptions are the notation for translation functors and the right-hand action of the affine Weyl group.

The roots of $B$ correspond to the negative roots of a root system $R\subset X(T)$, whose positive roots and simple roots are $R^{+}$ and $S$ respectively. 
The relations $\l-\alpha<\l$ for $\alpha\in R^{+}$ generate a partial ordering on $X(T)$.
$E$ is the Euclidean space $X(T)\otimes_{\mathbb Z}\mathbb R$  with its invariant bilinear form $\langle-, -\rangle$. $\l\in X(T)$ will be identified with $\l\otimes 1\in E$. $\rho$ is half the sum of the positive roots or, since $G$ is simply connected, the sum of the fundamental weights. The sets of dominant weights and restricted weights are respectively 
\begin{equation*}
\begin{aligned}
X(T)_{+}&=\{\lambda\in X(T): \langle\lambda,\,\alpha^{\scriptscriptstyle\vee}\rangle\ge0\hbox{ for all }\alpha\in R^{\,+}\}
\cr
X_r(T)&=\{\lambda\in X(T)_{+}: \langle\lambda,\,\alpha^{\scriptscriptstyle\vee}\rangle<p^r\hbox{ for all }\alpha\in S\}
\end{aligned}
\end{equation*}
where $\alpha^{\scriptscriptstyle\vee}$ is the dual root of $\alpha$. 

The connected components of the set $\{\lambda\in E: \langle\lambda+\rho,\,\alpha^{\scriptscriptstyle\vee}\rangle\not\in p\mathbb{Z}\hbox{ for all }\alpha\in R^{\,+}\}$ are called alcoves. Fix an alcove $C$ in the dominant region (not necessarily the standard alcove). $\overline{C}$ is bounded by a set of root hyperplanes. $\Sigma(C)$ is the set of reflections in these root hyperplanes. The affine Weyl group $W_p$ is the group generated by $\Sigma(C)$.  The natural action of $w\in W_p$ on $\l\in E$ is denoted $w\cdot\l$.

Suppose that $s\in\Sigma(C)$. The $s$-wall of $\overline{C}$ is the wall fixed by $s$. If $C'$ is any alcove then $C'=w_1\cdot C$ for a unique element $w_1\in W_p$. The $s$-wall of $\overline{C'}$ is the conjugate of the $s$-wall of $\overline{C}$ under the action of $w_1$. If $\l\in C'$ then $\l_s$ is the reflection of $\l$ in the $s$-wall of $\overline{C'}$. Note that $(w\cdot\l)_s=w\cdot(\l_s)$ for any $w\in W_p\,$. 

The translation functors are usually defined relative to the standard alcove. It will be convenient to redefine them relative to the alcove $C$. 

Fix $\l\in C\cap X(T)_{+}$ and let $s\in\Sigma(C)$. The condition $p\ge h$, where $h$ is the Coxeter number of $R$, is implicit in the choice of $\l$. It also permits the choice of $\mu\in X(T)$ in the $s$-wall of $\overline{C}$ so that the stabilizer of $\mu$ in $W_p$ is generated by $s$. $\op Ts$ is the translation functor from $\l$ to $\mu$ and $\op Rs$ is the translation functor from $\mu$ to $\l$. The wall-crossing functor $\op\T s$ is the composition $\op Rs\op Ts\,$. Since $\op Ts$ and $\op Rs$ are exact functors, $\op\T s$ is exact.

In order to define these explicitly, and to establish a correspondence with the standard notation, there is a unique $w\in W_p$ so that $w\cdot C=C_0$ where $C_0$ is the fundamental alcove\hfill\break \cite[II.6.2]{J}. The translation functors are $\op Ts=T_{w\cdot\l}^{w\cdot\mu}$ and $\op Rs=T_{w\cdot\mu}^{w\cdot\l}$ \cite[II.7.6]{J}. With this set-up, $\l_s=w^{-1}s_0w\cdot\l$, where $s_0$ is the reflection in the root hyperplane coincident with the $s$-wall of $C_0$. Consequently $\l_s<\l$ if, and only if, $w^{-1}s_0\cdot(w\cdot\l)<w^{-1}\cdot(w\cdot\l)$.

If $\nu\in X(T)_+$ then the simple module with highest weight $\nu$ is denoted by $L(\nu)$ and $H^0(\nu)$ denotes the global sections of a line bundle defined on $G/B$ with character $\nu$. 

\vfill\eject
\section{Translation-related Functors}
\label{tf}

Suppose that $V$ is a $G$ module. The image of $\epsilon_V$ is written $V_s$. The mappings $V\mapsto V_s$ on modules and $f\mapsto f_s$ on morphisms, where $f_s$ is the restriction of $f$ to $V_s$, define a functor. The properties of this functor and its relation to $\cop s{}$ were developed by the author in \cite{KC2} where the reader may find proofs of the results that are only stated in this section.
 
The functor $V\mapsto V_s$ is neither left exact nor right exact, though it does preserve injections and surjections \cite[2.8]{KC2}. The functor $\cop s{}$ is a right exact functor with right adjoint $\Ker\eta'$\hfill\break\cite[1.1]{KC2}. 

Standard properties of adjunction maps yield the following annihilation properties\hfill\break\cite[1.2]{KC2}.

\begin{prop}\indent\par
\label{prop: annihilate}
\begin{enumerate}
\item $\op Ts{}$ annihilates $\cop sV\,$.
\smallskip
\item $\cop s{}$ annihilates $\op RsV\,$.
\end{enumerate}
\end{prop}

If $L$ is a simple module, $\cop sL=0$ if, and only if, $\op TsL\ne0$. In turn, $\op TsL\ne0$ if, and only if, $L=L(w\cdot\l)$ where $w\in W_p$ and $w\cdot\l_s>w\cdot\l$ \cite[II.7.15]{J}. This leads to a vanishing criterion for $\cop s{}$ \cite[1.4]{KC2}.

\begin{prop}
\label{prop: cvanish}
$\cop sV=0$  if, and only if, $\cop sL=0$ for each simple quotient, $L\/$, of\/ $V$ .
\end{prop}

The following result characterizes $V_s$ as a submodule of $V$ \cite[1.8]{KC2}.

\begin{prop}\indent\par
\label{prop: char}
\begin{enumerate}
\item $\op Ts(V/V_s)=0$ and $V_s$ is the smallest submodule of $\,V$ with this property. 
\smallskip
\item $\cop sV_s=0$ and $V_s$ is the largest submodule of $\,V$ with this property. 
\end{enumerate}
\end{prop}

There is a unique $p$-regular weight $\l_0\in X_1(T)$ with $\l=\l_0+p\l_1$, $\l_1\in X(T)_+$. The next result of \cite{KC2} needed here relates $V\mapsto V_s$ to the analogously defined functor relative to $\l_0$.

The block of $G$ containing $L(\l)$ corresponds to the orbit of $\l$ under $W_p$. Let $\p_\l$ denote projection to the block of $L(\l)$.  Note that $L(\l)$ and $L(\l_0)$ may lie in different blocks. Suppose that $C'$ is the alcove containing $\l_0$. Choose $t\in\Sigma(C)$ so that the $s$-wall of $\overline{C}$ is the translate of the $t$-wall of $\overline{C'}$ by $p\l_1\,$.

\begin{prop}
\label{prop: transfer}
\textup{($p>h$)} Suppose that $V$ and $W$ are $G$ modules with $V$ finite-dimensional. Then
\[
((\p_{\l_0}V)\otimes W^{[1]})_s=\p_\l(V_t\otimes W^{[1]}).
\]
\end{prop}

It should be noted that the restriction on $p$ in \cite[3.5]{KC2} is stated differently here. There $G$ was simple and the condition $G\ne SL_p$ was sufficient to ensure that $p$ did not divide the root index $(X(T):\mathbb ZR)$. The condition $p>h$ is also sufficient and is better suited to the semisimple case.

Let $Q_\l$ denote the injective hull of $L(\l)$ as a $G$-module. Proposition \ref{prop: transfer} can be applied to $Q_\l$ when $p>2h-2$. This was used in \cite{KC2} to prove the following properties of $(Q_\l)_s$ \cite[3.1 and 3.8 where $I(\l)=Q_\l$]{KC2}.

\begin{prop}
\label{prop: Qs}
\textup{($p>2h-2$)}\par
\begin{enumerate}
\item $(Q_\l)_s$ is non-zero.
\smallskip
\item $(Q_\l)_s=Q_\l$ if, and only if, $\l_s>\l$.
\end{enumerate}
\end{prop}

\section{Indirectly Derived Functors}
\label{indirect}

The results of this section show how the right derived functors of $\kop s{}$ can be determined using cohomology of the functor $V\mapsto V_s\,$.
First, as was mentioned in the introduction, I need to show that $(\dop s,\kop s)$ is an adjoint pair of functors. 

If $V$ and $W$ are $G$ modules and $f\in\h_G(\op TsV,\op TsW)$ then $\alpha(f)=R_s(f)\,\eta'_V$ and $\beta(f)=\epsilon_W\,R_s(f)$ define adjunction isomorphisms \cite[A.6.2]{W}. If $\gamma=\beta\alpha^{-1}$ then $\gamma(\alpha(f))\,\eta'_V=\beta(f)\,\eta'_V=\epsilon_W\,\alpha(f)$. This leads to a commutative diagram,
\begin{equation}
\label{adj}
\begin{CD}
0@>>>\h_G(V,\kop sW)@>>>\h_G(V,\op\T sW)@>>>\h_G(V,W)\\
@.@.@V{\gamma}VV@|\\
0@>>>\h_G(\dop sV,W)@>>>\h_G(\op\T sV,W)@>>>\h_G(V,W).
\end{CD}
\end{equation}
Since the rows are exact, an isomorphism fills in the second column.

\begin{prop}
\label{prop: adjpair}
$\dop s$ is left adjoint to $\kop s$.
\end{prop}

By standard properties of adjoints \cite[2.6.1]{W}, $\kop s{}$ is left exact and $\dop s{}$ is right exact.

Since $\epsilon_V$ factors through $V_s$, there are two short exact sequences,
\begin{gather}
\label{def1}
0\longrightarrow\kop sV\longrightarrow\op\T sV\longrightarrow V_s\longrightarrow0\\
\label{def2}
0\longrightarrow V_s\longrightarrow V\longrightarrow \cop sV\longrightarrow0\,,
\end{gather}
which are natural in $V$.

Suppose that $(I\strut^q, d\strut^q)$ is a ($\kop s{}$-)acyclic resolution of a $G$-module $V$ with augmentation map $\iota\mathrel{:}V\rightarrow I^0$. Applying equation \ref{def1} to the resolution yields an exact sequence of complexes,
\[
0\longrightarrow\kop sI^{\textstyle\cdot}\longrightarrow\op\T sI^{\textstyle\cdot}\longrightarrow I^{\textstyle\cdot}_s\longrightarrow0\,.
\]
Because $\op\T s{}$ is exact, the long exact sequence of cohomology collapses in higher degrees yielding isomorphisms, $R\strut^q[\kop s{}]V\simeq H\strut^{q-1}(I_s^{\textstyle\cdot})$, for each $q>1$. In addition, the lower degree terms form an exact sequence,
\begin{equation*}
\begin{CD}
0@>>>H^0(\kop sI^{\textstyle\cdot})@>>>H^0(\op\T sI^{\textstyle\cdot})@>{\delta}>>H^0(I_s^{\textstyle\cdot})@>>>R^1[\kop s{}]V@>>>0\,.
\end{CD}
\end{equation*}

By naturality, there is a commutative diagram,
\begin{equation*}
\begin{CD}
\op\T sV@>{\epsilon_V}>>V_s\\
@V{\op\T s(\iota)}VV@VV{\iota_s}V\\
H^0(\op\T sI^{\textstyle\cdot})@>{\delta}>>H^0(I_s^{\textstyle\cdot})\,,
\end{CD}
\end{equation*}
where $\op\T s(\iota)$ is an isomorphism because $\op\T s{}$ is exact and $\iota_s$ is injective \cite[2.8]{KC2}.

Since $\delta$ factors through $V_s$, there is an exact sequence,
\begin{equation*}
\begin{CD}
0@>>>V_s@>>>H^0(I_s^{\textstyle\cdot})@>>>R^1[\kop s{}]V@>>>0\,.
\end{CD}
\end{equation*}
The following proposition summarizes these results.

\begin{prop}\indent\par
\label{prop: centralseq}
\begin{enumerate}
\item $R^1[\kop s{}]V\simeq H^0(I_s^{\textstyle\cdot})/V_s\,$.
\smallskip
\item $R\strut^q[\kop s{}]V\simeq H\strut^{q-1}(I_s^{\textstyle\cdot})$ for all $q>1\,$.
\end{enumerate}
\end{prop}

Applying \ref{def2} to the resolution yields another exact sequence of complexes,
\[
0\longrightarrow I_s^{\textstyle\cdot}\longrightarrow I^{\textstyle\cdot}\longrightarrow \cop sI^{\textstyle\cdot}\longrightarrow0\,.
\]
Here the long exact sequence begins with an injection of
$H^0(I_s^{\textstyle\cdot})$ into $V$. Recalling the definition of $\cop s{}$, there is an exact sequence,
\begin{equation}
\label{copseq}
0\longrightarrow R^1[\kop s{}]V\longrightarrow \cop sV\longrightarrow V/H^0(I_s^{\textstyle\cdot})\longrightarrow0\,.
\end{equation}

In particular, this proves the following result.

\begin{prop}
\label{prop: der1vanish}
If $\cop sV=0$ then $R^1[\kop s{}]V=0$.
\end{prop}

Proposition \ref{prop: centralseq} will now be used to derive some vanishing properties for the higher derived functors of~$\kop s{}$.

Because $\op Ts{}$ is an exact functor, proposition \ref{prop: char} implies that $\op TsV_s\simeq\op TsV$ and
\[
\op TsH\strut^q(I_s^{\textstyle\cdot})\simeq
H\strut^q(\op TsI_s^{\textstyle\cdot})\simeq
H\strut^q(\op TsI^{\textstyle\cdot})\simeq
\op TsH\strut^q(I^{\textstyle\cdot}).
\]
So $\op TsH^0(I_s^{\textstyle\cdot})\simeq\op TsV$ and $\op Ts{}$ applied to any higher cohomology vanishes.

On the other hand, by \ref{prop: annihilate}, $\cop s{}$ annihilates $\op RsV$ so that $(\op RsV)_s=\op RsV$. $\op RsI^{\textstyle\cdot}$ is an injective resolution of $\op RsV$ and
\[
H\strut^q((\op RsI^{\textstyle\cdot})_s)\simeq
H\strut^q(\op RsI^{\textstyle\cdot})\simeq
\op RsH\strut^q(I^{\textstyle\cdot}).
\]
Hence $H^0((\op RsI^{\textstyle\cdot})_s)\simeq\op RsV=(\op RsV)_s$ and any higher cohomology vanishes.

Applying \ref{prop: centralseq} yields the following result.

\begin{prop}\indent\par
\label{prop: vanish}
\begin{enumerate}
\item $\op Ts{}$ annihilates $R\strut^q[\kop s{}]V$ if $q>0$.
\smallskip
\item $\op RsV$ is acyclic.
\end{enumerate}
\end{prop}

It will be useful to have a condition that implies $R^1[\kop s{}]V\simeq\cop sV$. According to equation \ref{copseq}, this holds if, and only if, $H^0(I_s^{\textstyle\cdot})\simeq V$. If $\iota(V)\subseteq I_s^0$ then a diagram chase shows that there is an exact sequence,
\[
\begin{CD}
0@>>>V@>>>I_s^0@>d_s^0>>I_s^1\,.
\end{CD}
\]
Hence $V\simeq H^0(I_s^{\textstyle\cdot})\,$.

\begin{prop}
\label{prop: der1cop}
If $\iota(V)\subseteq I_s^0$ then $R^1[\kop s{}]V\simeq\cop sV$.
\end{prop}

A stronger condition will determine $R^2[\kop s{}]V$ as well. Let $B=I^0/\iota(V)$. Then, since $I^0$ is acyclic, $R^2[\kop s{}]V\simeq R^1[\kop s{}]B$. If $I_s^0=I^0$ then the condition of \ref{prop: der1cop} is satisfied and $\cop sB=0$, since $\cop s{}$ is right exact. By \ref{prop: der1vanish}, $R^1[\kop s{}]B=0$, which yields the following corollary.

\begin{cor}
\label{cor: der2cop}
If $I_s^0=I^0$ then
$R\strut^q[\kop s{}]V\simeq
\begin{cases}
\hfill\cop sV,\hfill&q=1\cr
\hfill0,\hfill&q=2\,.\cr
\end{cases}$
\end{cor}

This can be generalized to higher degrees. Let $B=I\strut^q/\im d\strut^{q-1}$ where $q\ge1$. Then the degree shift yields $R\strut^{q+2}[\kop s{}]V\simeq R^1[\kop s{}]B$. If $(I\strut^q)_s=I\strut^q$ then $R\strut^{q+2}[\kop s{}]V=0$.

\section{Acyclic Modules}
\label{inj}

The main aim of the section is to identify some $\kop s{}$-acyclic modules. First, we need to review the structure of $Q_\l$.

The injective hull of $L(\l_0)$ as a $G_1T$-module is $\widehat{Q}_1(\l_0)$. When $p\ge 2h-2$, $\widehat{Q}_1(\l_0)$ has a unique $G$-module structure and $Q_\l\simeq\widehat{Q}_1(\l_0)\otimes Q_{\l_1}^{[1]}\,$ \cite[II.11.16a]{J}.
Also, since $\widehat{Q}_1(\l_0)$ is indecomposable, it belongs to the block of $L(\l_0)$.

Assume that $p>2h-2$ for the remainder of this article. Note that this implies that $p>h$ so the restriction on the prime in \ref{prop: transfer} is satisfied. 

\begin{lem}
If $W$ is any $G$-module, then $\widehat Q_1(\l_0)\otimes W^{[1]}$ is acyclic for $\kop s{}$.
\end{lem}

If $W=0$, the result is trivial. Assuming that $W\ne0$, choose an injective resolution $0\rightarrow W\rightarrow Q^{\textstyle\cdot}$. The functor $V\mapsto V^{[1]}$ is exact so $(Q^{\textstyle\cdot})^{[1]}$ is an exact complex resolving $W^{[1]}$. Each tensor product $\widehat Q_1(\l_0)\otimes (Q\strut^q)^{[1]}$ is an injective G-module \cite[II.11.16a]{J}. Tensoring with $\widehat Q_1(\l_0)$  results in an injective resolution,
\begin{equation}
\label{resolution}
0\longrightarrow \widehat Q_1(\l_0)\otimes W^{[1]}\longrightarrow I^{\textstyle\cdot}\,,
\end{equation}
where $I^{\textstyle\cdot}=\widehat Q_1(\l_0)\otimes (Q^{\textstyle\cdot})^{[1]}$.

Using proposition \ref{prop: transfer} and the fact that $\widehat{Q}_1(\l_0)$ is finite-dimensional \cite[II.11.4]{J}, there is an equality of submodules,
\[
I_s^q=
\p_\l(\widehat Q_1(\l_0)_t\otimes(Q^q)^{[1]})\,,
\]
for each $q$. Since the coboundary maps are restrictions of the coboundary maps of $I^{\textstyle\cdot}$ to these submodules, this is also an equality of complexes. This shows that the cohomology of the complex $I^{\textstyle\cdot}_s$ vanishes in all non-zero degrees and that
\[
H^0(I^{\textstyle\cdot}_s)\simeq \p_\l(\widehat Q_1(\l_0)_t\otimes W^{[1]})=(\widehat Q_1(\l_0)\otimes W^{[1]})_s\,.
\]
Since an injective resolution is certainly an acyclic resolution, $\widehat Q_1(\l_0)\otimes W^{[1]}$ is acyclic by proposition \ref{prop: centralseq}. 

\begin{define}
\label{def: Adef}
Let $A_\l=\widehat Q_1(\l_0)\otimes H^0(\l_1)^{[1]}\,$.
\end{define}

Since $H^0(\l_1)$ is finite-dimensional and has simple socle $L(\l_1)$ \cite[I.5.12c, II.2]{J}, $H^0(\l_1)$ is isomorphic to a finite-dimensional submodule of $Q_{\l_1}$. So $A_\l$ is isomorphic to a finite-dimensional submodule of $Q_\l$. Using the lemma with $W=H^0(\l_1)$, $A_\l$ is acyclic. 

\begin{prop}
\label{prop: acycle}
$A_\l$ is acyclic.
\end{prop}

When $\l\in X_1(T)$, $\l=\l_0$ and $\l_1=0$ hence $A_\l=\widehat Q_1(\l)$.

\begin{cor}
\label{cor: Q1acycle}
If $\l\in X_1(T)$ then $\widehat Q_1(\l)$ is acyclic.
\end{cor}

If $\l\in X_r(T)$ and $r>1$, this generalizes to all $\widehat Q_r(\l)$ because $\widehat Q_r(\l)\simeq \widehat Q_1(\l_0)\otimes \widehat Q_{r-1}(\l_1)^{[1]}$ as $G_rT$ modules \cite[II.11.16b]{J}, and also as $G$ modules by the restriction on $p$.

\begin{cor}
\label{cor: Qracycle}
If $\l\in X_r(T)$ then $\widehat Q_r(\l)$ is acyclic.
\end{cor}

The properties of $(A_\l)_s$ that will be needed later are summarized in the following proposition (compare with \ref{prop: Qs}). 

\vfill\eject\begin{prop}\indent\par
\label{prop: Als}
\begin{enumerate}
\item $(A_\l)_s$ is non-zero.
\smallskip
\item $(A_\l)_s=A_\l$ if, and only if, $\l_s>\l$.
\end{enumerate}
\end{prop}

Recall that $\widehat Q_1(\l_0)$ belongs to the block of $L(\l_0)$ and, by the comments following \ref{def: Adef}, $A_\l$ belongs to the block of $L(\l)$. Applying \ref{prop: transfer}, $(A_\l)_s=\widehat Q_1(\l_0)_t\otimes H^0(\l_1)^{[1]}$. The proposition will follow from the corresponding properties of $\widehat Q_1(\l_0)$.

Let $\l_0^*=w_0\,\l_0+2(p-1)\rho$, which is the highest weight of $\widehat Q_1(\l_0)$ \cite[II.11.6a]{J}. The corresponding $H^0(\l_0^*)$ occurs at the top of the good filtration of $\widehat Q_1(\l_0)$ \cite[II.11.13]{J}.

Under the (twisted) duality $\tau$ defined in \cite[II.2.12]{J}, $\strut^{\tau}\widehat Q_1(\l_0)$ is the projective cover of $L(\l_0)$ as a $G_rT$ module. By \cite[II.11.5.3]{J}, $\widehat Q_1(\l_0)$ is self-dual and so the head of $\widehat Q_1(\l_0)$ is $L(\l_0)$. Since the Weyl module $V(\l_0^*)\simeq\strut^{\tau}H^0(\l_0^*)$ \cite[II.2.13.2]{J}, $\widehat Q_1(\l_0)$ has a submodule isomorphic to $V(\l_0^*)$.

There are two cases to consider.

If $\l_s>\l$ then $(\l_0)_t>\l_0$. By proposition \ref{prop: cvanish}, $\cop tL(\l_0)=0$ implies that $\cop t\widehat Q_1(\l_0)=0$. Then $\widehat Q_1(\l_0)_t=\widehat Q_1(\l_0)$ by proposition \ref{prop: char} and $(A_\l)_s=A_\l$.

If $\l_s<\l$ then $(\l_0)_t<\l_0$. Since $\op TtL(\l_0)=0$, \ref{prop: char} implies that $\widehat Q_1(\l_0)_t\ne \widehat Q_1(\l_0)$ which completes the proof of part (ii). Also $(\l_0^*)_t>\l_0^*$ so that $\op TtL(\l_0^*)\ne 0$ implies that $\cop tV(\l_0^*)=0$ by \ref{prop: cvanish}. Using \ref{prop: char}, $V(\l_0^*)$ injects into $\widehat Q_1(\l_0)_t$. By \ref{prop: transfer}, $V(\l_0^*)\otimes H^0(\l_1)^{[1]}$ injects into~$(A_\l)_s\,$. 

In both cases, $(A_\l)_s$ is non-zero which proves (i).

As an example, consider the (partial) calculation of $R\strut^q[\kop s{}]L(\l)$. $\kop s L(\l)=0$ if $\l_s<\l$, because $\op TsL(\l)=0$. When $\l_s>\l$, $\kop sL(\l)=\rad\op\T sL(\l)$ \cite[II.7.19]{J}. In either case, because $(A_\l)_s\subseteq Q_\l$ is non-zero, there is an isomorphism
\[
\h_G(L(\l),(A_\l)_s)\simeq\h_G(L(\l),A_\l)\simeq k\,.
\]
Any augmentation map $\iota: L(\l)\rightarrow A_\l$ must factor through $(A_\l)_s$ so the condition of \ref{prop: der1cop} is satisfied. Then $R^1[\kop s{}]L(\l)\simeq \cop sL(\l)$. If $\l_s<\l$, $\cop sL(\l)\simeq L(\l)$. If $\l_s>\l$ then $\cop sL(\l)=0$ and $(A_\l)_s=A_\l$. Applying \ref{prop: der1cop} and \ref{cor: der2cop} gives the following results.

\begin{ex}
\label{ex: Lex}
\indent\par
\begin{enumerate}
\item If $\l_s<\l$, 
$R\strut^q[\kop s{}]L(\l)\simeq
\begin{cases}
\hfill0,\hfill&q=0\cr
\hfill L(\l),\hfill&q=1\,.\cr
\end{cases}$
\medskip
\item If $\l_s>\l$, 
$R\strut^q[\kop s{}]L(\l)\simeq
\begin{cases}
\hfill\rad\op\T sL(\l),\hfill&q=0\cr
\hfill0,\hfill&q=1\text{ or\, }2\,.\cr
\end{cases}$
\end{enumerate}
\end{ex}

These results imply that \ref{prop1} and \ref{prop2} are not true in general. If $\l_s<\l$ then $\cop sA_\l$ is non-zero by proposition \ref{prop: Als} while $A_\l$ is acyclic. Let $V=A_\l/L(\l)$. Because $A_\l$ is acyclic, $R^2[\kop s{}]V\simeq R^1[\kop s{}]L(\l)$. If $\l_s<\l$, then $R^2[\kop s{}]V\simeq L(\l)\ne0$.

\section{Line Bundle Calculations}
\label{lb}

The calculation of the derived functors $R\strut^q[\kop s{}]H^0(\l)$ will be the subject of this section. The first objective is to show that the hypothesis of \ref{prop: der1cop} holds when $V=H^0(\l)$.

\begin{lem}
\label{lemma: sinclusion}
$\h_G(H^0(\lambda),V(\l_0^*)\otimes H^0(\l_1)^{[1]})\simeq k$.
\end{lem}

By standard properties of injective envelopes, the dimension of $\h_G(H^0(\l),Q_\l)$ is the multiplicity of $L(\l)$ in $H^0(\l)$. Hence $\h_G(H^0(\l),Q_\l)$ is one-dimensional. Any non-zero subspace is isomorphic to $k$ and, since $H^0(\l)$ injects into $Q(\l)$, any non-zero homomorphism is injective. To prove the lemma, it will suffice to show that the homomorphism space is non-zero.

Using the fact that $V(\l_0^*)^*\simeq H^0(-w_0^{\vphantom*}\l_0^*)$,
\begin{equation*}
\h_G(H^0(\l),V(\l_0^*)\otimes H^0(\l_1)^{[1]})\simeq\h_G(H^0(\l)\otimes H^0(-\l_0 +2(p-1)\rho),H^0(\l_1)^{[1]}).
\end{equation*}

The tensor product $H^0(\l)\otimes H^0(-\l_0 +2(p-1)\rho)$ has a good filtration and, since the highest weight is $2(p-1)\rho+p\l_1$, its corresponding $H^0$ is at the top of the filtration. Now it suffices to show that $\h_G(H^0(2(p-1)\rho+p\l_1),H^0(\l_1)^{[1]})$ is non-zero.

This is a special case of a result Jantzen used in the proof of \cite[II.11.13.4]{J}.
There he shows, for any $\nu\in X(T)_+$, that there is a non-zero $G$-homomorphism from $H^0(\nu+(p-1)\rho)$ to $L(\nu')\otimes H^0(\nu_1)^{[1]}$ where $\nu'=w_0\nu_0+(p-1)\rho$ . If $\nu=(p-1)\rho+p\l_1$ then  $\nu+(p-1)\rho=2(p-1)\rho+p\l_1$ and $\nu'=0$ which completes the proof of the lemma.

Since $V(\l_0^*)\otimes H^0(\l_1)^{[1]}$ is isomorphic to a submodule of $A_\l$, the lemma implies that $\h_G(H^0(\lambda),A_\l)\ne0$.

\begin{cor}
\label{cor: bottom}
$\h_G(H^0(\lambda),A_\l)\simeq k$.
\end{cor}

Since any two non-zero homomorphisms are injective and are scalar multiples, $A_\l$ has a unique submodule isomorphic to $H^0(\l)$. Also,
an acyclic resolution of $H^0(\l)$ can be constructed starting with any non-zero homomorphism $\iota: H^0(\l)\rightarrow A_\l$. Showing that such an $\iota$ satisfies the hypothesis of \ref{prop: der1cop} comes down to proving the following proposition.

\begin{prop}
\label{prop: sinclusion}
$\h_G(H^0(\lambda),(A_\l)_s)\simeq \h_G(H^0(\lambda), A_\l)$.
\end{prop}

If $\l_s>\l$, $(A_\l)_s=A_\l$ so there is nothing to prove. When $\l_s<\l$, $V(\l_0^*)\otimes H^0(\l_1)^{[1]}$ injects into~$(A_\l)_s\,$. By the lemma, $\h_G(H^0(\lambda),(A_\l)_s)$ is non-zero.

Since the homomorphism spaces in the proposition are one-dimensional, the image of $\iota$ is contained in $(A_\l)_s$ which is the condition of \ref{prop: der1cop}.

\begin{prop}
\label{prop: R1}
$R^1[\kop s{}]H^0(\l)\simeq\cop sH^0(\l)$.
\end{prop}

This result will be used to determine $R^1[\kop s{}]H^0(\l)$. The calculation of $\kop s{}H^0(\l)$ will be developed at the same time. The results depend on the position of the $s$-wall of $\overline{C}$, that is, whether $\l_s<\l$ or $\l_s>\l$.

If $\l_s<\l$ there are two possibilities.

If $\l_s\in X(T)_{+}$ then \cite[Lemma 3.2]{A} yields a non-split exact sequence,
\begin{equation*}
0\longrightarrow H^0(\l_s)\longrightarrow \op RsH^0(\mu)\longrightarrow H^0(\l)\longrightarrow 0\,,
\end{equation*}
where $\op Rs{}$ is the translation functor from $\mu$ to $\l$\, as defined in \S1.

Because $\op TsH^0(\l)\simeq H^0(\mu)$ \cite[2.1a]{A}, this can be re-written
\begin{equation}
\label{ses}
0\longrightarrow H^0(\l_s)\longrightarrow \op\T sH^0(\l)\longrightarrow H^0(\l)\longrightarrow 0\,.
\end{equation}
Since $\epsilon_{H^0(\l)}$ is non-zero and $\h_G(\op\T sH^0(\l), H^0(\l))\simeq\h_G(H^0(\mu),H^0(\mu))\simeq k$ \cite[II.2.8]{J}, this implies that $\kop s{}H^0(\l)\simeq H^0(\l_s)$ and $H^0(\l)_s=H^0(\l)$ hence $R^1[\kop s{}]H^0(\l)=0$.

If $\l_s$ is not in $X(T)_{+}$, $\op TsH^0(\l)=0$ so that $\kop s{}H^0(\l)=0$ and $H^0(\l)_s=0$. Then $R^1[\kop s{}]H^0(\l)\simeq H^0(\l)$. 
Note that, in this situation, $\l_s=s_\alpha\cdot\l$, where $\alpha$ is a simple root and $0<\langle \l+\rho,\,\alpha^{\scriptscriptstyle\vee}\rangle<p$.
By \cite[II.5.4d]{J}, $H^0(\l)\simeq H^1(\l_s)$.

The two results can be combined into a single statement since $H^0(\l_s)=0$ if $\l_s$ is not in $X(T)_{+}$ \cite[II.2.6]{J} and $H^1(\l_s)=0$ for $\l_s\in X(T)_{+}$ by Kempf's vanishing theorem \cite[II.4.5]{J}.

\begin{thm}
\label{thm: sRqI}
If $\l_s<\l$ then $R\strut^q[\kop s{}]H^0(\l)\simeq H\strut^q(\l_s)$ for $q\le1$.
\end{thm}

It seems likely that the theorem is true for all $q\ge0$. (See conjecture \ref{conj: der2vanish} below.)

If $\l_s>\l$, then \cite[Theorem 3.3]{A} implies that there is a non-zero intertwining homomorphism $\Psi\mathrel{:}H^0(\l_s)\rightarrow H^0(\l)$ which is unique up to scalar multiple. Let $K(\l_s)=\Ker\Psi$ and $C(\l)=\Coker\Psi$. 
Because $\op TsH^0(\l)\simeq\op TsH^0(\l_s)$ and $\op Ts{}$ is exact, $\op TsK(\l_s)=0$ and $\op TsC(\l)=0$. Since $\cop s{}$ is right exact and $\cop sH^0(\l_s)=0$, 
$\cop s{}(\im\Psi)=0\,$. Hence $H^0(\l)_s=\im\Psi$ by \ref{prop: char}. Proposition \ref{prop: R1} then yields the following result.

\begin{thm} 
\label{thm: R1s}
If $\l_s>\l$ then $R^1[\kop s{}]H^0(\l)\simeq C(\l)$.
\end{thm}

The calculation of $\kop sH^0(\l)$ in this case requires a little more work. There is a commutative diagram where the rows are exact and the vertical maps are natural.
\begin{equation}
\label{R0s}
\begin{CD}
0@>>>\op\T sK(\l_s)@>>>\op\T sH^0(\l_s)@>>>\op\T sH^0(\l)_s@>>>0\cr
@.@VVV@VVV@VVV@.\cr
0@>>>K(\l_s)@>>>H^0(\l_s)@>>>H^0(\l)_s@>>>0\cr
\end{CD}
\end{equation}
Note that $K(\l_s)\ne0$, since otherwise $\soc_G H^0(\l)$ is not simple, and $\op\T sK(\l_s)=0$ because $\op TsK(\l_s)=0$. An application of the snake lemma yields the exact sequence,
\begin{equation}
\label{non-split}
0\longrightarrow H^0(\l)\longrightarrow \kop sH^0(\l)\longrightarrow K(\l_s)\longrightarrow 0\,,
\end{equation}
because $\kop sH^0(\l_s)\simeq H^0(\l)$ by \ref{thm: sRqI} applied to $\l_s$. Since $\soc_G\op\T sH^0(\l)$ is simple \cite[II.7.19b]{J}, $\kop sH^0(\l)$ is indecomposable and so this sequence is non-split. To show that $\kop sH^0(\l)$ is the unique non-split extension up to isomorphism it suffices to show that the dimension of $\E_G^1(K(\l_s),H^0(\l))$ is at most one.

Using the adjoint pairing $(\op Ts{},\op Rs{})$, $\E_G^q(K(\l_s), \op\T sH^0(\l))=0$ for all $q\ge0$ (since $\op TsK(\l_s)=0$). Interchanging the roles of $\l$ and $\l_s$ in \ref{ses} results in an exact sequence,
\begin{equation}
\label{sess}
0\longrightarrow H^0(\l)\longrightarrow \op\T sH^0(\l_s)\longrightarrow H^0(\l_s)\longrightarrow0\,.
\end{equation}
Applying $\h_G(K(\l_s),-)$ to \ref{sess} yields an isomorphism,
\[
\h_G(K(\l_s),H^0(\l_s))\simeq\E_G^1(K(\l_s),H^0(\l)). 
\]
Using Frobenius reciprocity \cite[I.3.4]{J},
\[
\dim\h_G(K(\l_s),H^0(\l_s))=\dim\h_B(K(\l_s),\l_s)\le\dim\h_T(K(\l_s),\l_s).
\]
Because $K(\l_s)$ is a submodule of $H^0(\l_s)$, $\dim\E_G^1(K(\l_s),H^0(\l))\le1$.

\begin{thm}
\label{thm: R0s}
If $\l_s>\l$ then $\kop sH^0(\l)$ is the unique non-split extension of $K(\l_s)$ by $H^0(\l)$ up to isomorphism.
\end{thm}

Taken together the  three previous results give the values of $R\strut^q[\kop s{}]H^0(\l)$ for $q\le1$. The next result relates the higher derived functors of the two cases. By \ref{prop: vanish}ii, $\op\T sH^0(\l)$ is acyclic. The long exact sequence obtained by applying $\kop s{}$ to  equation \ref{sess} collapses in higher degrees proving the following.

\begin{prop}
\label{prop: increment}
If $\l_s>\l$ then $R\strut^{q+1}[\kop s{}]H^0(\l)\simeq R\strut^q[\kop s{}]H^0(\l_s)$ for all $q\ge1$.
\end{prop}

In particular, since $R^1[\kop s{}]H^0(\l_s)\simeq H^1(\l)=0$, the proposition leads to the following corollary. 

\begin{cor}
\label{cor: R2}
If $\l_s>\l$ then $R^2[\kop s{}]H^0(\l)=0$.
\end{cor}

The results of this section and the parallel results in category $\mathcal O$, particularly \ref{prop2}, suggest the following conjecture.  

\begin{conj}
\label{conj: der2vanish}
$R\strut^q[\kop s{}]H^0(\l)=0$ for all $q>1$.
\end{conj}

The conjecture has two equivalent formulations.

\begin{prop}The following statements are equivalent:
\begin{enumerate}
\item $R\strut^q[\kop s{}]H^0(\l)=0$ for all $q>1$.
\smallskip
\item $\displaystyle A_\l/H^0(\l)$ is acyclic.
\smallskip
\item $R\strut^q[\kop s{}]H^0(\l)\simeq H\strut^q(\l_s)$ for all $q\ge0$ if $\l_s<\l$.
\end{enumerate}
\end{prop}

Since $A_\l$ is acyclic, $R\strut^q[\kop s{}](A_\l/H^0(\l))$ is isomorphic to $R\strut^{q+1}[\kop s{}]H^0(\l)$ for $q\ge1$. Hence $A_\l/H^0(\l)$ is acyclic if, and only if, conjecture \ref{conj: der2vanish} holds.

Conjecture \ref{conj: der2vanish} immediately implies that the statement of theorem \ref{thm: sRqI} is true for all $q\ge0$ since $H\strut^q(\l_s)=0$ when $q\ge2$. Conversely, suppose that theorem \ref{thm: sRqI} is true for all $q\ge0$. When $\l_s<\l$, $R\strut^q[\kop s{}]H^0(\l)=0$ for $q\ge2$. Suppose that $\l_s>\l$. Applying \ref{prop: increment} together with \ref{cor: R2}, $R\strut^{q+1}[\kop s{}]H^0(\l)=0$ for $q\ge1$.

\section{Extensions}
\label{ext}

In this section, the preceding ideas are applied to extensions of $G$ modules. The main goal is to derive some consequences of conjecture \ref{conj: der2vanish}. A secondary goal is to consider the existence of a Grothendieck spectral sequence \cite[5.8.3]{W} corresponding to the composition $\h_G(V,\kop s(-))$ where $\eta'_V:V\rightarrow \op\T sV$ is an injective map (compare with \cite[2.10ii]{KC1}). The sufficient condition for the existence of this spectral sequence is that $\E_G^q(V,\kop sQ)=0$ for all $q>0$ and for all injective $Q$. The first step in this direction is the following proposition.

\begin{prop}
\label{prop: E1}
Suppose that $\eta'_V$ is an injective map and that $Q$ is an injective module.
\begin{enumerate}
\item$\h_G(V,Q_s)\simeq\h_G(V,Q)$.
\smallskip
\item$\E_G^1(V,\kop sQ)=0$.
\smallskip
\item$\E_G^q(V,\kop sQ)\simeq\E_G^{q-1}(V,Q_s)$ for all $q>1$.
\end{enumerate}
\end{prop}

Applying \ref{adj} with $W=Q$ yields a commutative diagram,
\begin{equation}
\label{Qadj}
\begin{CD}
0@>>>\h_G(V,\kop sQ)@>>>\h_G(V,\op\T sQ)@>>>\h_G(V,Q)\\
@.@VVV@V{\gamma}VV@|\\
0@>>>\h_G(\dop sV,Q)@>>>\h_G(\op\T sV,Q)@>>>\h_G(V,Q).
\end{CD}
\end{equation}

\vfill
\eject
The last horizontal map in the first row of \ref{Qadj} is $\h_G(V,\epsilon_Q)$. Let $i:Q_s\rightarrow Q$ be the natural inclusion.  Since $\epsilon_Q$ factors through $Q_s$, $\epsilon_Q=i\,(\epsilon_Q)_s$ and there is a commutative diagram,
\begin{equation}
\label{fulladj}
\begin{CD}
0@>>>\h_G(V,\kop sQ)@>>>\h_G(V,\op\T sQ)@>{\beta}>>\h_G(V,Q_s)@>>>0\\
@.@VVV@V{\gamma}VV@VV{\alpha}V@.\\
0@>>>\h_G(\dop sV,Q)@>>>\h_G(\op\T sV,Q)@>>>\h_G(V,Q)@>>>0\,,
\end{CD}
\end{equation}
where $\alpha=\h_G(V,i)$ and $\beta=\h_G(V,(\epsilon_Q)_s)$.

The second row of \ref{fulladj} is exact because $\eta'_V$ is an injective map and $Q$ is an injective module. Because $\gamma$ is an isomorphism, the composition $\alpha\,\beta$ is surjective. This means that $\alpha$ is surjective. Since $\alpha$ is injective, $\alpha$ is an isomorphism (proving (i)). Note that $\beta$ is surjective so the first row of this diagram is also exact.

The mapping $\beta$ fits into a long exact sequence,
\[
\begin{CD}
0@>>>\h_G(V,\kop sQ)@>>>\h_G(V,\op\T sQ)@>{\beta}>>\h_G(V,Q_s)\\
@.@>>>\E_G^1(V,\kop sQ)@>>>\E_G^1(V,\op\T sQ)@>>>\cdots
\end{CD}
\]
Since $\beta$ is surjective and $\op\T sQ$ is injective, $\E_G^1(V,\kop sQ)=0$ (proving (ii)). The same long exact sequence yields isomorphisms, $\E_G^{q-1}(V,Q_s)\simeq\E_G^q(V,\kop sQ)$, for all $q>1$ (proving (iii)).

The next result is an important special case that will be used to draw some consequences from conjecture \ref{conj: der2vanish}. It may not be the best possible such result (see \ref{prop: SS} below). The first part of the proof is based on ideas from \cite[5.1.8]{GJ}. 

\begin{thm}
\label{thm: SS0}
Suppose that $\eta'_V$ is injective and $W$ is acyclic for $\kop s{}$. Then there is a commutative diagram, for all $q\ge1$,
\[
\begin{CD}
\cdots@>>>\E_G^{q-1}(V,W)@>>>\E_G^q(\dop sV,W)@>>>\E_G^q(\op\T sV,W)@>>>\cdots\\
@.@V{\alpha_{q-1}}VV@V{\gamma_q}VV@V{\beta_q}VV@.\\
\cdots@>>>\E_G^{q-1}(V,W_s)@>>>\E_G^q(V,\kop sW)@>>>\E_G^q(V,\op\T sW)@>>>\cdots
\end{CD}
\]
where the $\beta_q$ are isomorphisms and the $\alpha_q$ and $\gamma_q$ are injective maps. If, in addition, $W_s$ is acyclic then the $\alpha_q$ and the $\gamma_q$ are isomorphisms.
\end{thm}

Suppose that $I^{\textstyle\cdot}$ is an injective resolution of $W$. Applying \ref{def1} to the resolution yields an exact sequence of cochain complexes,
\begin{equation*}
0\longrightarrow\kop sI^{\textstyle\cdot}\longrightarrow\op\T sI^{\textstyle\cdot}\longrightarrow I^{\textstyle\cdot}_s\longrightarrow0\,.
\end{equation*}
Since $W$ is acyclic, \ref{prop: centralseq} implies that these are each resolutions of the terms in the exact sequence,
\begin{equation}
\label{Wses}
0\longrightarrow\kop sW\longrightarrow\op\T sW\longrightarrow W_s\longrightarrow0\,.
\end{equation}
Let $W_s\rightarrow J^{\textstyle\cdot}$ and $\kop sW\rightarrow L^{\textstyle\cdot}$ be injective resolutions. Using the horseshoe lemma \cite[2.2.8 (in the opposite category)]{W}, there is an injective resolution $\op\T sW\rightarrow K^{\textstyle\cdot}$ which completes a split exact sequence resolving \ref{Wses},
\begin{equation}
\label{horseshoe}
0\longrightarrow L^{\textstyle\cdot}\longrightarrow K^{\textstyle\cdot}\longrightarrow J^{\textstyle\cdot}\longrightarrow0\,.
\end{equation}
Since $I_s^{\textstyle\cdot}$ and $\kop sI^{\textstyle\cdot}$ are resolutions, the comparison theorem \cite[2.3.7]{W} yields cochain maps $a:I_s^{\textstyle\cdot}\rightarrow J^{\textstyle\cdot}$ and $c:\kop sI^{\textstyle\cdot}\rightarrow L^{\textstyle\cdot}$. Using the injectivity of $L^{\textstyle\cdot}$ and the splitting maps of \ref{horseshoe}, there is a cochain map $b:\op\T sI^{\textstyle\cdot}\rightarrow K^{\textstyle\cdot}$ completing a commutative diagram of cochain complexes with exact rows,
\begin{equation*}
\begin{CD}
0@>>>\kop sI^{\textstyle\cdot}@>>>\op\T sI^{\textstyle\cdot}@>>>I^{\textstyle\cdot}_s@>>>0\\
@.@V{c}VV@V{b}VV@V{a}VV@.\\
0@>>>L^{\textstyle\cdot}@>>>K^{\textstyle\cdot}@>>>J^{\textstyle\cdot}@>>>0\,.
\end{CD}
\end{equation*}
Applying $\h_G(V,-)$ yields
\begin{equation*}
\begin{CD}
0@>>>\h_G(V,\kop sI^{\textstyle\cdot})@>>>\h_G(V,\op\T sI^{\textstyle\cdot})@>>>\h_G(V,I^\cdot_s)@.\\
@.@VVV@VVV@VVV@.\\
0@>>>\h_G(V,L^{\textstyle\cdot})@>>>\h_G(V,K^{\textstyle\cdot})@>>>\h_G(V,J^{\textstyle\cdot})@>>>0\,.
\end{CD}
\end{equation*}

According to \ref{fulladj}, there is also a commutative diagram with exact rows,
\begin{equation*}
\begin{CD}
0@>>>\h_G(V,\kop sI^{\textstyle\cdot})@>>>\h_G(V,\op\T sI^{\textstyle\cdot})@>>>\h_G(V,I^\cdot_s)@>>>0\\
@.@VVV@VVV@VV{\h_G(V,\,i)}V@.\\
0@>>>\h_G(\dop sV,I^{\textstyle\cdot})@>>>\h_G(\op\T sV,I^{\textstyle\cdot})@>>>\h_G(V,I^{\textstyle\cdot})@>>>0\,,
\end{CD}
\end{equation*}
where $i: I^{\textstyle\cdot}_s \rightarrow I^{\textstyle\cdot}$ is the natural inclusion and the vertical maps are isomorphisms. Combining the last two diagrams results in a commutative diagram,
\begin{equation*}
\begin{CD}
0@>>>\h_G(\dop sV,I^{\textstyle\cdot})@>>>\h_G(\op\T sV,I^{\textstyle\cdot})@>>>\h_G(V,I^{\textstyle\cdot})@>>>0\\
@.@V{\gamma}VV@V{\beta}VV@V{\alpha}VV@.\\
0@>>>\h_G(V,L^{\textstyle\cdot})@>>>\h_G(V,K^{\textstyle\cdot})@>>>\h_G(V,J^{\textstyle\cdot})@>>>0\,.
\end{CD}
\end{equation*}
The cochain version of \cite[1.3.4]{W} establishes the first claim. (It will be convenient to use the same letter for the induced maps on cohomology.)

Since $\op\T sW\rightarrow\op\T sI^{\textstyle\cdot}$ is an injective resolution, $b$ is a homotopy equivalence. Then $\h_G(V,\,b)$ induces an isomorphism on cohomology and so $\beta$ will be an isomorphism.

Let $f: J^{\textstyle\cdot}\rightarrow I^{\textstyle\cdot}$ be a cochain map lifting the inclusion $W_s\rightarrow W$. The composition $fa$ is a cochain map from $I_s^{\textstyle\cdot}$ to $I^{\textstyle\cdot}$ lifting the inclusion $W_s\rightarrow W$. Then, by the uniqueness part of the comparison theorem, $fa$ is homotopic to $i$ . Since $\h_G(V,i)$ is an isomorphism of complexes, $\h_G(V,a)$ induces an injection on cohomology which means that $\alpha$ is also injective.

The mapping $\gamma_0$ is an isomorphism by \ref{prop: adjpair}, so it is certainly injective. For each $q\ge 1$, $\alpha_{q-1}$ and $\beta_q$ are injective and $\beta_{q-1}$ is surjective. The long five lemma \cite[1.3.3]{W} implies that each $\gamma_q$ is injective.

Finally, if $W_s$ is acyclic then $J_s^{\textstyle\cdot}$ is a resolution of $W_s$. Arguing as above, $af_s$ is homotopic to the natural inclusion $j: J_s^{\textstyle\cdot}\rightarrow J^{\textstyle\cdot}$. Since $\h_G(V,j)$ is an isomorphism, $\h_G(V,a)$ induces a surjection on cohomology. Then $\h_G(V,a)$ induces an isomorphism on cohomology and $\alpha$ is an isomorphism. The long five lemma establishes the last claim.

If conjecture \ref{conj: der2vanish} is true then theorem \ref{thm: SS0} will produce some consequences: 

Suppose that $w\cdot\l\in X(T)_+$ for some $w\in W_p$. Recall that $\l\in X(T)_+$ was chosen in a fixed, but arbitrary, dominant alcove $C$. The translation functors were defined in such a way that any results expressed in terms of $\l$ are equally true when expressed in terms of $w\cdot\l$. If $w\cdot\l_s>w\cdot\l$ then, by \ref{thm: sRqI}, $\kop sH^0(w\cdot\l_s)=H^0(w\cdot\l)$ and $R^1[\kop s{}]H^0(w\cdot\l_s)=0$.  Conjecture \ref{conj: der2vanish} would imply that  $H^0(w\cdot\l_s)$ is acyclic so that \ref{thm: SS0} could be applied with  $W=H^0(w\cdot\l_s)$. By 
\ref{prop: R1}, $\cop sH^0(w\cdot\l_s)=0$ which means that $W_s=W$ would be acyclic.

If $\l_s>\l$ then there are two substitutions that can be made for $V$ in \ref{thm: SS0}. First, let $V=H^0(\l)$. Since $\op\T sH^0(\l)\simeq\op\T sH^0(\l_s)$, \ref{ses} can be re-written as
\begin{equation*}
0\longrightarrow H^0(\l)\longrightarrow \op\T sH^0(\l)\longrightarrow H^0(\l_s)\longrightarrow 0\,.
\end{equation*}
Since $\eta'_V$ is non-zero and $\h_G(H^0(\l), \op\T sH^0(\l))\simeq\h_G(H^0(\mu),H^0(\mu))\simeq k$ \cite[II.2.8]{J}, $\eta'_V$ is injective and $\dop sV\simeq H^0(\l_s)$. Second, let $V=L(\l)$. Since $\op\T sV$ is non-zero, $\eta'_V$ is non-zero so must also be injective and $\dop sV\simeq\op\T sL(\l)/L(\l)$.

\begin{prop}
\label{prop: consequences}
Assume that conjecture \ref{conj: der2vanish} is true.
Let $w\cdot\l\in X(T)_+$ for some $w\in W_p$. Suppose that $\l_s>\l$ and $w\cdot\l_s>w\cdot\l$. Then, for all $q\ge0\,$,
\begin{enumerate}
\item $\E_G^q(H^0(\l),H^0(w\cdot\l))\simeq\E_G^q(H^0(\l_s),H^0(w\cdot\l_s))$\text{ and }
\smallskip
\item $\E_G^q(L(\l),H^0(w\cdot\l))\simeq\E_G^q(\op\T sL(\l)/L(\l),H^0(w\cdot\l_s))$.
\end{enumerate}
\end{prop}

The conclusions are reasonable and are consistent with known results on extensions.  From that point of view, they lend some indirect support to conjecture \ref{conj: der2vanish}.

Further consequences of \ref{conj: der2vanish} along these lines could be derived if the spectral sequence mentioned above exists. At this point, the existence question remains open. We conclude by giving a sufficient condition for the existence expressed as a conjecture.

\begin{conj}
\label{conj: Qsacyclic}If $w\in W_p$ and $w\cdot\l\in X(T)_+$ then
$(Q_{w\cdot\l})_s$ is acyclic for $\kop s{}$.
\end{conj}

Note that, by \ref{prop: Qs}, $(Q_{w\cdot\l})_s=Q_{w\cdot\l}$ when ${w\cdot\l}_s>{w\cdot\l}$. Since $Q_{w\cdot\l}$ is injective, it is also acyclic. So, it would suffice to prove the conjecture when ${w\cdot\l}_s<{w\cdot\l}$. 

\begin{prop}
\label{prop: SS}Assume that conjecture \ref{conj: Qsacyclic} is true.
If $\eta'_V$ is an injective map then there is a spectral sequence,
\[
E_2^{r,q}=\E_G^r(V,R\strut^q[\kop s{}]W)\implies \E_G^{r+q}(\dop sV,W)\,.
\]
\end{prop}

If $Q$ is any injective module then $Q_s\subseteq\p_\l Q$ hence $Q_s=(\p_\l Q)_s$. Then $Q_s$ is a direct sum of modules of the form $(Q_{w\cdot\l})_s$ for $w\in W_p$ where $w\cdot\l\in X(T)_+$. Hence $Q_s$ is acyclic. Theorem \ref{thm: SS0} can be applied to $Q$ (since it is acyclic), so $\E_G^q(V,\kop sQ)\simeq\E_G^q(\dop sV,Q)=0$ for all $q>0$.  The composition $\h_G(V,\kop sW)$ is naturally isomorphic to $\h_G(\dop sV,W)$ by \ref{prop: adjpair} which completes the proof.

Note that, when $W$ is acyclic, the spectral sequence above collapses to an isomorphism,\hfill\break$\E_G^q(V,\kop sW)\simeq\E_G^q(\dop sV,W)$. This would be an improvement over the conclusion of \ref{thm: SS0} and would more directly imply the conclusions of \ref{prop: consequences}. So the existence of the spectral sequence would make conjecture \ref{conj: der2vanish} even more believable. 

It is not surprising that acyclicity of modules seems to be at the root of these problems. It is the underlying general concept that made the analogous results work in category $\mathcal O$ (although it was much easier going there). It is hoped that novel techniques will be found that will prove these conjectures for $G$ modules.

\vfill\eject

\end{document}